\documentclass[12pt]{article}

\usepackage[utf8]{inputenc}       
\usepackage[T1]{fontenc}          
\usepackage{lmodern}              
\usepackage{geometry}             
\usepackage{graphicx}             
\usepackage{amsmath, amssymb}     
\usepackage[numbers]{natbib}      
\usepackage{hyperref}             
\usepackage{caption}              
\usepackage{pgfplots}
\usepackage{tikz}
\usepackage{algorithm}
\usepackage{algpseudocode}

\newcommand{\epsrel}{\varepsilon_{\text{rel}}}
\newcommand{\epsabs}{\varepsilon_{\text{abs}}}
\newcommand{\mutarget}{\mu_\text{target}}

\title{Hybridizing PDHG and Interior-Point Methods}
\author{Edward Rothberg \\
Gurobi Optimization, LLC \\
\texttt{rothberg@gurobi.com}}
\date{March 3, 2026}

\begin{document}

\maketitle
\begin{abstract}
  The Primal-Dual Hybrid Gradient (PDHG) algorithm is a first-order
  method that can exploit GPUs to solve large-scale linear programming
  problems.  The approach can often be faster than the alternatives,
  simplex and interior-point methods, typically at the cost of much
  lower accuracy.  This paper looks at whether PDHG can be hybridized
  with an interior-point method to retain some of the speed advantages
  of the former while capturing the accuracy advantages of the latter.
\end{abstract}

\section{Introduction}

The Primal-Dual Hybrid Gradient (PDHG) algorithm~\cite{pdhg11} has
received a great deal of attention recently, primarily because it
provides an algorithmic framework for building a GPU-friendly method
to solve large linear programming (LP) problems.  Recent enhancements
to the method~\cite{pdlp22,pdhg11,hprlp24,cupdlp23,cupdlp+25} have made
it quite competitive with the current state-of-the-art algorithms for
solving such problems, the simplex method~\cite{simplex51} and
interior-point methods~\cite{wright97}.  As a first-order method,
though, PDHG convergence can be quite slow.  As a result, it is
typically terminated at a point when the resulting solutions are
much less accurate than those computed by the alternatives.

This paper looks at whether PDHG can be hybridized with an
interior-point method to obtain some of the advantages of both
methods.  Specifically, we compute a low-accuracy solution using PDHG
and use that solution to warm-start an interior-point method.
Interior-point methods are notoriously difficult to warm
start~\cite{warmstart2010,warmstart1998,warmstart2008,wright97,warmstart2002},
but this reputation was formed in a different context.  Earlier work
tried to use an optimal solution for one LP to accelerate the solution
of a similar LP.  The task of starting from an
approximate solution to the same problem has not received much
attention, presumably because methods like PDHG that can produce
near-optimal solutions quickly have not been competitive before.

There are of course other options for obtaining more accuracy from a
PDHG solution.  The most obvious is to tighten the PDHG convergence
criteria, giving it more time to improve the accuracy of its solution.
Another is to feed the PDHG solution into a
crossover~\cite{crossover91} procedure to perform a simplex-like
procedure that computes a much more accurate {\em basic\/} solution.
Each presents its own tradeoffs.  Runtime versus accuracy is one, of
course; another is the degree to which each is amenable to multi-core
or GPU parallelism.  This paper will consider several such
alternatives.

We find that an inexpensive strategy for starting an interior-point
algorithm from a low-accuracy PDHG solution is both quick and
reliable.  The average reduction in iteration counts, when
compared with a cold-started method, was roughly 3X, with the
warm-started version solving in fewer than 10 iterations more than
half of the time.  The resulting interior-point solutions have 3 to 4
orders of magnitude smaller violations than the PDHG solutions they
start from.

\section{Background}
\label{sec:Background}

A linear programming problem can be stated in its primal and dual
forms as:
\begin{equation*}
\begin{minipage}{0.45\linewidth}
\begin{align*}
\min \quad & c^\top x \\
\text{s.t.} \quad & A x = b \\
& x \ge 0
\end{align*}
\end{minipage}
\hspace{0.2cm}
\begin{minipage}{0.45\linewidth}
\begin{align*}
\max \quad & b^\top y \\
\text{s.t.} \quad & A^\top y + z = c \\
                  & z \ge 0
\end{align*}
\end{minipage}
\end{equation*}
Constraint matrix $A$ has $m$ rows and $n$ columns, objective vector
$c$ has length $n$, and right-hand-side vector $b$ has length $m$.
Vectors $x$ and $y$ are the primal and dual solutions, and vector $z$
provides {\em reduced costs\/} ($z_j$ gives the rate of change on the
objective per unit change in $x_j$).

Practical linear programming models are typically more complex than
this canonical form, often involving inequalities and non-trivial
lower and upper bounds on variables.  These can be easily mapped into
the form above (and are typically handled implicitly in LP solver
implementations).

Any trio of vectors $(x,y,z)$ has associated primal and dual residual
vectors:
\begin{align*}
r_P &= b - Ax \\
r_D &= c - A^\top y - z
\end{align*}
A solution is optimal if (i) the primal and dual residual vectors are
zero, (ii) $x$ and $z$ are non-negative, and (iii) the objective gap is
zero. This last condition can be expressed directly as
$c^\top x = b^\top y$, or somewhat less directly as $x^\top z = 0$
(referred to as {\em complementarity\/}).

The PDHG algorithm maintains a current iterate, and
iterations proceed until all three optimality
conditions are met to target
tolerances (typically expressed in relative terms):
\begin{align*}
  \| r_P \|_2 & \leq \epsrel (1 + \| b \|_2) \\
  \| r_D \|_2 & \leq \epsrel (1 + \| c \|_2) \\
  | c^\top x - b^\top y| & \leq \epsrel (1 + |c^\top x| + |b^\top y|)
\end{align*}
The most commonly used convergence tolerance when evaluating the
performance of these methods is $\epsrel = 10^{-4}$.

PDHG is a first-order method, which in practice means that convergence
can be quite slow.  The strongest claim that has been made for the
approach is that convergence is linear when the problem exhibits
certain properties.  Tighter convergence tolerances typically lead to
significant increases in runtime.

Interior-point methods also maintain a current iterate that does not
necessarily satisfy feasibility or complementarity targets.
Interior-point methods are {\em locally quadratically convergent\/},
meaning that they converge quite quickly once the iterate is ``close
enough'' to an optimal solution.  Convergence theory puts strong
bounds on the number of iterations required to reach
convergence~\cite{wright97}.
These bounds are based on relative tolerances, but practical
implementations typically iterate well beyond the tolerances that
convergence theory would suggest.

Tolerances for simplex and crossover methods are typically expressed
in absolute terms ($\|r_P\|_{\infty} \leq \epsabs$ and
$\|r_D\|_{\infty} \leq \epsabs$), with the traditional target being
$\epsabs = 10^{-6}$.  The (primal) simplex algorithm maintains
(primal) feasibility and zero complementarity (to the extent that this
is numerically possible), so its iterations focus solely on achieving
one of these three optimality conditions.  When compared with PDHG,
the smaller absolute target and the use of infinity norms means that
simplex termination conditions are much stricter.

\subsection{Presolve and Scaling}
\label{sec:preprocessing}

In a practical LP solver, the model solved by the LP algorithm is
typically quite different from the model provided by the user.  The
two main differences are due to {\em presolve\/}~\cite{presolve1995}
and {\em scaling\/}~\cite{scaling2011,scaling1975}.

Presolve is a process of identifying properties of the variables and
constraints that are true in any solution, and exploiting them to
create an equivalent model that is smaller and easier to solve.  For
example, if it can be shown that $x_1 = 1$ or $x_1 = x_2$ in an
optimal solution, then $x_1$ can be removed from the model before it
is passed to the solver.  Presolve is an important part of any modern
LP solver, and it typically reduces the cost of solving the model
substantially.

Scaling is a process of multiplying each row and column in the
constraint matrix by a constant to improve numerical conditioning.
The process of solving an LP typically involves linear algebra on
large matrices, and better conditioning leads to smaller errors in
these computations, which typically improves speed and robustness.
Scaling is helpful for interior-point methods, important for the
simplex method, and essential for PDHG.

The reason we mention these problem transformations is to point out
that the constraint violations in the solution found by the LP solver
will not be the same as those in the solution returned to the user.
While these transformations could in theory make the errors smaller or
larger, unfortunately they almost always make them larger, often
substantially so.  As a result, the need for the solver to produce an
accurate solution is not driven entirely by the user, but also by the
fact that accuracy can degrade substantially when a solution is
translated back to the user model.

\section{Warm-Starting an Interior-Point Method}

Our hope with a method that hybridizes PDHG and interior-point methods
is to take advantage of the strengths of both: have PDHG compute a
quick, low-accuracy solution, and create from that a near-optimal
starting point for the interior-point solver that can take advantage
of its quadratic convergence near the optimal solution.

On the surface, the steps in this approach are straightforward:
perform PDHG iterations on the model until a (loose) target tolerance
is achieved, start the interior-point solver from that solution, and
then run the interior-point solver until convergence.

The difficulty here lies in the fact that an interior-point solver
cannot be warm-started from just any point.  Convergence depends
heavily on how well {\em centered\/} the starting point is.  An
interior-point solver actually solves a series of optimization
problems, where the goal in each is to achieve primal feasibility,
dual feasibility, and $x_j \cdot z_j = \mu$ for all $j$ (for some
$\mu$) ~\cite{wright97}.  The algorithm drives $\mu \rightarrow 0$ as
it proceeds, pushing the complementarity term $x^\top z$ to 0 and thus
pushing the solution towards an optimal point (assuming primal and
dual feasibility have also been achieved).  In practical
implementations, $\mu$ is reduced after each iteration of the method,
so the problem associated with a fixed $\mu$ is never fully solved.
However, interior-point convergence proofs rely on step length
guarantees that only hold if the iterate remains sufficiently close to
a point where $x_j \cdot z_j = \mu$ for all $j$ (the set of such
points is called the {\em central path\/}).

Unfortunately, an arbitrary start point is unlikely to be primal
feasible, dual feasible, and also meet this centering condition.  For
example, the start point may be an optimal solution to a related model
that has been modified (e.g., after adding a variable bound or cutting
plane when solving a Mixed-Integer Programming
problem~\cite{warmstart1998}), so the solution is no longer feasible.
Or variables in the solution may be exactly at their bounds, giving
$x_j \cdot z_j = 0$ (as is common in a PDHG solution).  Even when
neither are true, having $x_j \cdot z_j \approx \mu$ for some $\mu$ is
a fairly esoteric condition that is unlikely to be met by an
arbitrary start vector.  Yildirim and Wright~\cite{warmstart2002} lay
out the conditions under which a start vector will produce a favorable
outcome, and the conditions are quite strict.

Despite these challenges, careful construction of a suitable starting
point has shown some promise in practical implementations.
Gondzio~\cite{warmstart1998} demonstrated single-digit iteration
counts within a cutting plane application.  John and
Yildirim~\cite{warmstart2008} found improvements in iteration counts,
but these typically did not translate into reduced runtimes due to the
costs of computing an effective starting point.  Engau, Anjos, and
Vannelli demonstrated a 50\% reduction in iteration counts over a
fairly broad test set.  To our knowledge, though, no widely used
interior-point LP solver implements warm-starting, so the techniques
developed so far have not had a big impact on the practice of solving
LPs.

\subsection{Our Warm-Start Approach}

While the literature describes several sophisticated (and sometimes
expensive) approaches for computing effective warm-start vectors, as
we noted earlier, we suspect that the problem we are trying to solve is
simpler than the one considered before: previous approaches were
specifically trying to reoptimize in the face of non-trivial changes
to the underlying problem, whereas we are trying to start from a
near-optimal solution to the unmodified problem.

We tried several different schemes for constructing a warm-start vector
from a PDHG solution, and eventually settled on a simple and
inexpensive approach, described in Algorithm~1.
\begin{algorithm}
\caption{Centered start point}
\begin{algorithmic}[1]
    \State \textbf{Input:} PDHG solution $(x, y, z)$ to scaled LP, $\mu_{\min}$, $\delta_{\max}$, $\alpha_{\min}$
    \State $\mutarget \gets \max(x^\top z / n, \mu_{\min})$
    \ForAll{$j \in \{1..n\}$}
      \State $x'_j \gets \max(x_j, \alpha_{\min})$, $z'_j \gets \max(z_j, \alpha_{\min})$
      \If{$x'_j < z'_j$}
        \State $x'_j \gets \min(\max(\mutarget / z'_j,\,x'_j - \delta_{\max}),\,x'_j + \delta_{\max})$
        \State $z'_j \gets \min(\max(\mutarget / x'_j,\,z'_j - \delta_{\max}),\,z'_j + \delta_{\max})$
      \Else
        \State $z'_j \gets \min(\max(\mutarget / x'_j,\,z'_j - \delta_{\max}),\,z'_j + \delta_{\max})$
        \State $x'_j \gets \min(\max(\mutarget / z'_j,\,x'_j - \delta_{\max}),\,x'_j + \delta_{\max})$
      \EndIf
    \EndFor
    \State Warm-start the interior-point solver from $(x', y, z')$.
    \While{interior-point step is very small}
      \State $\alpha_{\min} \gets 10 \cdot \alpha_{\min}$
      \ForAll{$j \in \{1..n\}$}
        \State $x'_j \gets \max(x'_j, \alpha_{\min})$, $z'_j \gets \max(z'_j, \alpha_{\min})$
      \EndFor
    \EndWhile
\end{algorithmic}
\end{algorithm}
We experimented with multiple values for the input parameters
and found that setting the minimum complementarity
target $\mu_{\min}$ to $10^{-6}$, the minimum variable value
$\alpha_{\min}$ to $10^{-6}$, and
the maximum variable perturbation
$\delta_{\max}$ to $10^{-4}$ provided the best results.

The details of the approach were definitely important.  We found that
the most important one was to work with the scaled model.  One goal of
scaling is to create a constraint matrix $A$ where the norms of each
row and column are roughly equal to 1.  One consequence of this property
is that a constant perturbation to a variable produces a constant
increase in the primal and dual infeasibility (in infinity norms).
Large increases in infeasibility make it difficult to achieve quick
termination.

Another important detail was the need for a strategy to avoid
getting stuck in the first interior-point iteration.  We
simultaneously want our starting point to be as close as possible to
the PDHG solution to accelerate convergence, but also not so close to
a boundary that the solver is unable to take a substantial step.  We
found that increasing $\alpha_{\min}$ until the interior-point solver
was able to take one non-trivial step struck a nice balance between
these two goals.  Requiring a minimum step length of $10^{-6}$ worked
well in our tests.

One modification we tried to this centering scheme seemed promising
but was ultimately unsuccessful.  Once structural variables were
perturbed, we tried recomputing primal and dual slacks.  These derived
values were then perturbed to improve centering.  While we did find
that this approach produced smaller residuals, convergence was
significantly worse.

\section{Experimental Results}

We have a number of options for solving linear programs, each with its
own tradeoffs between runtime and accuracy.  We will present results for
a number of those options here.  Before doing so, we first describe
our testing environment, including details on the software, machines,
and LP test sets that were used.

\subsection{Testing Environment}

All of our tests were performed using Gurobi version 13.0, extended to
implement the methods described here.  The PDHG implementation in this
version of Gurobi closely follows the approach used in
cuPDLP+~\cite{cupdlp+25}.  Gurobi also includes state-of-the-art
implementations of the interior-point and crossover algorithms, which
have been tested and tuned over many years and a wide range of
practical LP models.

Most of our tests were performed on an Nvidia Grace Hopper GH200
system, which consists of an Nvidia Grace CPU containing 72 cores and
a Hopper H100 GPU.  In our tests, PDHG runs on the GPU, and everything
else runs on the CPU.  The system has 512GB of shared memory, and the
GPU has its own 96GB of HBM memory.  This is a unique resource for us,
which limits our ability to do extensive GPU testing on a broad range
of large models.

To partly address this constraint, some of our tests were performed on
a large cluster of AMD EPYC 4364P-based systems, each with an 8-core
processor and 128~GB of memory.  As CPU-only systems with limited
memory bandwidth, these machines give quite poor performance for
PDHG.  We use them here simply because the large cluster
allows us to run tests on a much larger test set.

\subsection{Test Models}

The primary test set we use in this paper is the Mittelmann LP
set~\cite{hans25}, which consists of 43 publicly available models of
modest size and difficulty.  The mean solve time for these models is
only a few seconds.  All tests on this set were run on our GPU system
(PDHG on the GPU; everything else on the CPU).

We also consider a broader test set of 2298 LP models that capture a
mix of publicly available and customer models.  This broader set was
only run on our large CPU-based cluster.

We use a time limit of 10,000 seconds for all of our tests.  All
runtime results are measured as the time from passing the user model
to the solver to the time an optimal solution is returned.  When
presenting data on runtimes and solution accuracy, we sometimes
provide values for each individual model in the set, and sometimes
aggregate over the whole set.  Aggregation is done using geometric
means.

\subsection{Violation Metric}

One other form of aggregation we do in our testing is to combine
violations of primal feasibility, dual feasibility, and
complementarity into a single metric.  More precisely, the three
ingredients we use are the maximum primal infeasibility (measured as
$\|r_P\|_{\infty}$), the maximum dual infeasibility (measured as
$\|r_D\|_{\infty}$), and the objective gap (measured as
$x^\top z/(1 + |c^\top x|)$).  All are computed on the original model,
after undoing the effects of presolve and scaling.  We take the
maximum over these three quantities and call it the {\em maximum
  violation\/} for that solution to that model.

Why did we choose this violation measure?  Infinity norms are the
traditional metric used for capturing primal and dual infeasibility in
LP.  This matches with our own experience, where users typically
consider how much the solution violates any one constraint.  When it
comes to the objective gap, our experience is that users typically
consider the number of digits that agree.

We should add that we have performed these same comparisons using
other metrics, and while the violations change in absolute terms, the
relative results between methods are qualitatively similar.

\subsection{The Performance-Accuracy Tradeoff for Simple Variants}

Let us begin by comparing runtimes and maximum violations for three
readily-available options for improving on the solution accuracy
provided by a baseline PDHG run (with $\epsrel = 10^{-4}$).  The first
two involve simply tightening tolerances; in our case we
consider $\epsrel = 10^{-6}$ and $\epsrel = 10^{-8}$.  The third is to
use an interior-point method, without crossover.

\begin{figure}[htbp]
\centering
\begin{tikzpicture}
\begin{loglogaxis}[
    width=0.95\linewidth,
    height=0.62\linewidth,
    grid=both,
    legend style={at={(0.98,0.98)},anchor=north east},
    mark size=2.2pt,
    xmin=1,
    xmax=1e2,
    ymin=1e-12,
    ymax=1e6,
    xlabel={Relative runtime},
    ylabel={Max violation},
    ]

\definecolor{OIblue}{RGB}{0,114,178}
\definecolor{OIorange}{RGB}{230,159,0}
\definecolor{OIteal}{RGB}{0,158,115}
\definecolor{OIvermillion}{RGB}{213,94,0}

\addplot[
    only marks,
    mark=*,
    color=OIblue,
] table {
x y
1.000e+00 6.242e-04
1.000e+00 5.047e-05
1.021e+00 1.389e-03
1.000e+02 6.533e-03
9.983e+01 5.626e-05
9.870e+01 1.543e-02
1.000e+00 8.679e-03
1.000e+00 1.937e-01
1.000e+00 1.197e-03
1.000e+00 1.480e+00
1.000e+00 2.776e+00
1.000e+00 2.226e-05
1.716e+00 1.000e+06
1.000e+00 3.265e-05
1.000e+00 1.500e-03
1.000e+00 1.903e-03
1.979e+00 5.524e-01
1.000e+00 7.046e+00
1.000e+00 3.179e+01
1.062e+01 9.945e+02
1.000e+00 3.699e-05
1.000e+00 9.838e+01
1.104e+00 7.419e-01
1.619e+00 1.433e-01
1.000e+00 2.107e-03
1.000e+00 2.539e-03
1.000e+00 3.501e+00
1.070e+01 5.526e-03
2.293e+00 2.757e-02
4.940e+00 1.866e-01
1.000e+00 3.862e-04
1.000e+00 6.777e-04
1.306e+00 3.407e-01
4.257e+00 1.090e-02
1.000e+00 4.878e-01
1.000e+00 1.413e+00
1.000e+00 3.409e-03
1.000e+00 9.952e-04
1.000e+00 1.377e-02
1.000e+00 1.183e-02
3.252e+00 5.866e-03
1.000e+00 2.018e-02
1.000e+00 1.000e+06
};
\addlegendentry{PDHG $\epsrel = 10^{-4}$}

\addplot[
    only marks,
    mark=square*,
    color=OIorange,
] table {
x y
1.013e+00 2.895e-05
1.010e+00 9.326e-06
1.064e+00 1.109e-05
1.000e+02 1.265e-05
1.000e+02 4.004e-06
1.000e+02 2.198e-04
1.062e+00 2.844e-05
4.461e+00 8.322e-04
1.081e+00 9.800e-06
1.020e+00 9.382e-07
1.195e+00 2.337e-01
1.000e+00 1.000e-12
1.246e+01 1.000e+06
1.003e+00 3.368e-07
1.785e+00 6.671e-05
1.067e+00 1.466e-04
4.734e+00 5.410e-03
1.118e+00 3.216e-03
1.000e+02 6.202e-02
1.000e+00 4.975e-07
2.509e+00 1.415e+00
4.926e+00 2.287e-03
3.238e+00 4.129e-04
4.190e+00 3.930e-05
3.756e+00 2.125e-05
1.802e+00 2.490e-02
1.400e+01 8.704e-04
4.519e+00 1.864e-04
7.972e+00 4.589e-02
1.102e+00 5.293e-06
1.750e+00 7.997e-06
2.446e+00 7.117e-04
2.218e+01 3.418e-05
1.176e+01 1.744e-03
1.871e+00 6.863e-03
4.855e+00 2.173e-05
2.048e+00 6.510e-06
5.077e+00 1.627e-04
1.018e+00 1.586e-04
4.010e+00 5.910e-03
2.718e+00 1.050e-04
5.141e+00 4.759e+04
};
\addlegendentry{PDHG $\epsrel = 10^{-6}$}

\addplot[
    only marks,
    mark=triangle*,
    color=OIteal,
] table {
x y
1.016e+00 2.895e-05
1.038e+00 8.264e-07
1.085e+00 1.064e-07
1.000e+02 6.464e-06
1.000e+02 3.657e-07
1.000e+02 4.014e-06
1.109e+00 1.715e-06
1.387e+01 7.959e-04
1.135e+00 3.557e-07
1.023e+00 4.751e-05
1.732e+00 2.503e-04
1.000e+00 1.000e-12
2.323e+01 1.000e+06
1.009e+00 4.253e-08
2.252e+00 8.479e-07
1.067e+00 1.466e-04
5.432e+00 1.897e-04
1.145e+00 5.447e-05
1.000e+02 5.908e-04
1.000e+00 3.176e-09
2.625e+00 1.546e-02
6.025e+00 1.526e-05
6.667e+00 1.749e-06
5.267e+00 4.633e-07
1.914e+01 9.998e-07
2.970e+00 2.439e-04
2.535e+01 1.625e-05
7.387e+00 1.183e-05
1.047e+01 1.015e-02
1.119e+00 6.184e-08
4.013e+00 9.006e-07
3.911e+00 4.639e-06
2.917e+01 6.866e-07
1.000e+02 3.614e-05
3.073e+00 2.373e-05
1.087e+01 1.506e-07
2.270e+00 1.128e-06
1.209e+01 2.891e-06
1.036e+00 2.750e-06
5.267e+00 6.169e-07
9.904e+00 3.245e-06
9.838e+00 3.451e+05
};
\addlegendentry{PDHG $\epsrel = 10^{-8}$}

\addplot[
    only marks,
    mark=diamond*,
    color=OIvermillion,
] table {
x y
1.207e+00 1.077e-10
1.223e+00 1.933e-09
1.000e+00 5.103e-09
1.000e+00 3.817e-09
1.000e+00 3.690e-09
1.000e+00 4.904e-08
1.594e+00 3.945e-11
3.885e+00 9.512e-07
1.019e+01 1.428e-08
1.477e+00 4.366e-08
1.366e+00 2.116e-08
1.059e+00 1.412e-08
1.000e+00 1.000e+06
1.762e+00 3.126e-12
2.822e+00 2.139e-10
2.267e+00 1.000e-12
1.000e+00 1.084e-08
1.329e+00 1.177e-11
4.560e+00 6.482e-07
1.000e+00 9.153e-06
4.857e+00 1.678e-10
8.018e+00 2.296e-09
1.000e+00 1.758e-07
1.000e+00 9.226e-09
5.229e+00 1.744e-07
1.474e+00 9.202e-09
9.020e+00 9.356e-10
1.000e+00 3.138e-11
1.000e+00 2.070e-10
1.000e+00 1.460e-07
7.932e+00 1.610e-07
3.289e+00 2.856e-09
1.000e+00 2.359e-12
1.000e+00 1.885e-06
1.370e+01 2.244e-04
5.700e+00 4.669e-10
4.655e+00 7.922e-10
3.111e+00 3.915e-09
3.719e+00 6.723e-08
2.509e+00 5.741e-09
1.000e+00 3.248e-11
1.077e+01 4.160e-08
};
\addlegendentry{Interior point}
\end{loglogaxis}
\end{tikzpicture}
\caption{Runtimes and maximum violations for the models in the Mittelmann test set, for four solution approaches.  Runtime ratios larger than 100 are reported as 100.  Similarly, maximum violations smaller than $10^{-12}$ or larger than $10^{6}$ are reported as those values.}
\label{fig:count_vs_violation_base}
\end{figure}
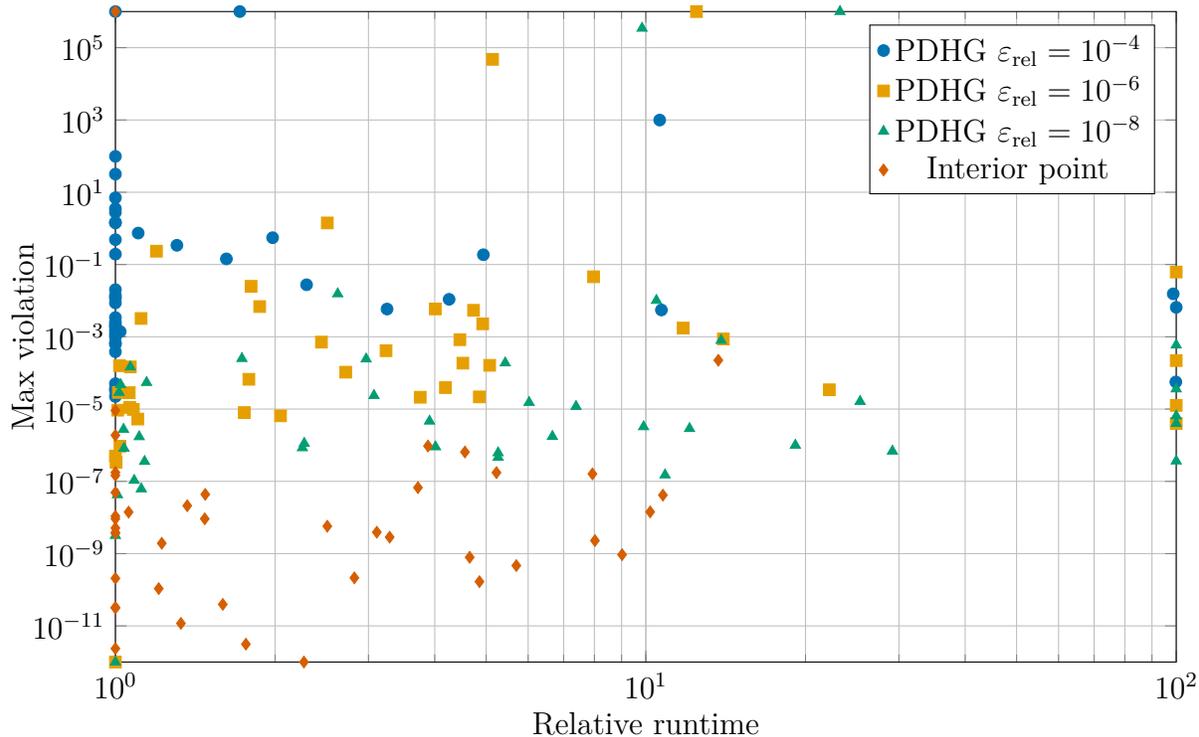

Figure~\ref{fig:count_vs_violation_base} shows the results of this
test on the Mittelmann LP test set (run on our Grace Hopper system).
Each point captures the results of a single method run on a single
model.  The x-axis shows the runtime for that model, relative to the
best runtime of any of the four methods on the same model.  The y-axis
shows the maximum constraint violation for that model (as defined
previously).


One thing that is clear from the results in this figure is that
tightening PDHG convergence tolerances has its limits;
violations are reduced, but at a substantial cost.  And even with very tight
tolerances, violations are still much larger than those for the
interior-pont solver.  Recall that interior-point methods are
{\em locally quadratically convergent\/}, so once the iterate gets close to
the optimal solution, each subsequent iterate brings with it a significant
improvement in solution accuracy.

Table~\ref{tab:mittelmann_basicmethods_aggregated} shows aggregated
results for this same test.
\begin{table}[htbp]
\centering
\begin{tabular}{c|c|c|c|c}
  \textbf{} &
   \shortstack{\textbf{PDHG}\\ $\epsrel = 10^{-4}$} &
   \shortstack{\textbf{PDHG}\\ $\epsrel = 10^{-6}$} &
   \shortstack{\textbf{PDHG}\\ $\epsrel = 10^{-8}$} &
   \shortstack{\textbf{Interior point}}\\ \hline
Models solved & 43 & 42 & 42 & 42 \\
Relative runtime & 1.00 & 2.29 & 3.56 & 1.16 \\
Mean max violation & $5 \cdot 10^{-2}$ & $3 \cdot 10^{-4}$ & $1 \cdot 10^{-5}$ & $1 \cdot 10^{-8}$\\
\end{tabular}
\caption{Runtime and maximum violation comparison, Mittelmann test: PDHG with different (relative) convergence tolerances, and interior point without crossover (computed using geometric means over the 43 models in the Mittelmann LP set).}
\label{tab:mittelmann_basicmethods_aggregated}
\end{table}
The first row shows the number of models that each method was able to
solve.  Only PDHG with $\epsrel = 10^{-4}$ was able to solve all 43
models to the desired tolerances.  The other three options are unable
to achieve target tolerances on one numerically challenging model
({\em ns1688926\/}).

The next row shows geometric means of total runtimes,
relative to the runtime for the baseline option (PDHG with
$\epsrel = 10^{-4}$).
The second shows geometric means of maximum constraint violations
for each approach.

A few clear conclusions can be drawn from these results: (i) PDHG
performance degrades significantly as tolerances are tightened, and
(ii) the interior-point solver produces violations that are 3-6 orders
of magnitude smaller than those produced by PDHG.  While it is
difficult to answer the question of how much accuracy is truly
necessary, the sheer magnitude of the difference suggests that some of
these methods will not provide sufficient accuracy for some
applications.

\subsection{Results for the Hybrid Approach}

We now consider the performance of our hybrid approach.  To set the
results in a broader context, we consider a number of additional
options for obtaining a solution that is more accurate than what you would get
from baseline PDHG (with $\epsrel = 10^{-4}$):
\begin{itemize}
\item Tighter PDHG tolerances; use $\epsrel = 10^{-8}$
  instead.  This option was also considered in the previous comparison.
\item Perform a crossover step, starting from the PDHG
  solution. Crossover runtimes can depend quite heavily on the
  accuracy of the starting solution, and we found that
  $\epsrel = 10^{-6}$ provided the best overall runtimes.
\item Warm-start an interior-point solver from the PDHG solution,
  using the strategy we just described.  We found
  that the best results were obtained by starting from
  a PDHG solution with $\epsrel = 10^{-4}$.
\end{itemize}



\begin{table}[t]
\centering
\begin{tabular}{c|c|c|c|c}
  \textbf{} &
   \shortstack{\textbf{PDHG}\\ \textbf{$\epsrel = 10^{-8}$}} &
   \shortstack{\textbf{PDHG}\\ $\epsrel = 10^{-6}$\\ + \textbf{crossover}} &
   \shortstack{\textbf{PDHG}\\ $\epsrel = 10^{-4}$\\ + \textbf{IP warm-start}} &
   \shortstack{\textbf{Interior point}}\\ \hline
Models solved & 42 & 42 & 35 & 42\\
Relative runtime & 3.56 & 2.84 & 2.17 & 1.16\\
Mean max violation & $1 \cdot 10^{-5}$ & $1 \cdot 10^{-10}$ & $2 \cdot 10^{-8}$ & $1 \cdot 10^{-8}$\\
\end{tabular}
\caption{Runtime and maximum violation comparison: PDHG with tighter (relative) convergence tolerances, PDHG with crossover, our hybrid approach, and interior point.  All are computed using geometric means over the 43 models in the Mittelmann LP test set.}
\label{tab:accuracy_options}
\end{table}
Results for these options for the Mittelmann LP test set,
run on our Grace Hopper system are summarized
in Table~\ref{tab:accuracy_options}.  We also
include results for a cold-started interior-point solver.
The first line shows the number
of models out of the 43 in the test set that were solved to the
requested tolerances with each method.  As noted earlier, most
approaches have trouble solving model {\em ns1688926\/}.  Our proposed
hybrid approach has trouble with a few more, solving 35 of the 43
models in the set.

As discussed earlier, the primary failure mode for our warm-started
approach is for it to get stuck, unable to take a substantial step
from the current iterate.  Our approach for adjusting the start point
helped a lot, but did not address the issue entirely.

The next line shows geometric-mean runtimes for the various methods,
relative to the baseline PDHG run (not shown), for the models that
each method was able to solve.  The PDHG-based approaches are much
more expensive than the baseline.  Interestingly, performing a
crossover step is faster than tightening PDHG tolerances, despite the
fact that crossover is a highly sequential operation, running on one
core of the CPU.  Our proposed hybrid approach is the fastest of the
alternatives, although it does not have the robustness of the others.

The final line shows geometric means of maximum violations.  This is
our first test that includes crossover, and the accuracy results are
substantiall better than those of any of the other options we
considered.  This is expected - the solutions produced by simplex
and crossover are known to be extremely accurate.  Note also that our
PDHG/interior-point hybrid significantly improves accuracy, as
expected.

One thing that is apparent from the chart is that the cold-started
interior-point solver is the overall winner on this set.  We will
return to this topic shortly.

We want to reiterate that different applications have different
accuracy requirements, and that more accuracy is not necessarily
better.  Our intent is just to point out that different methods can
provide radically different solution qualities, which could have
significant downstream implications.


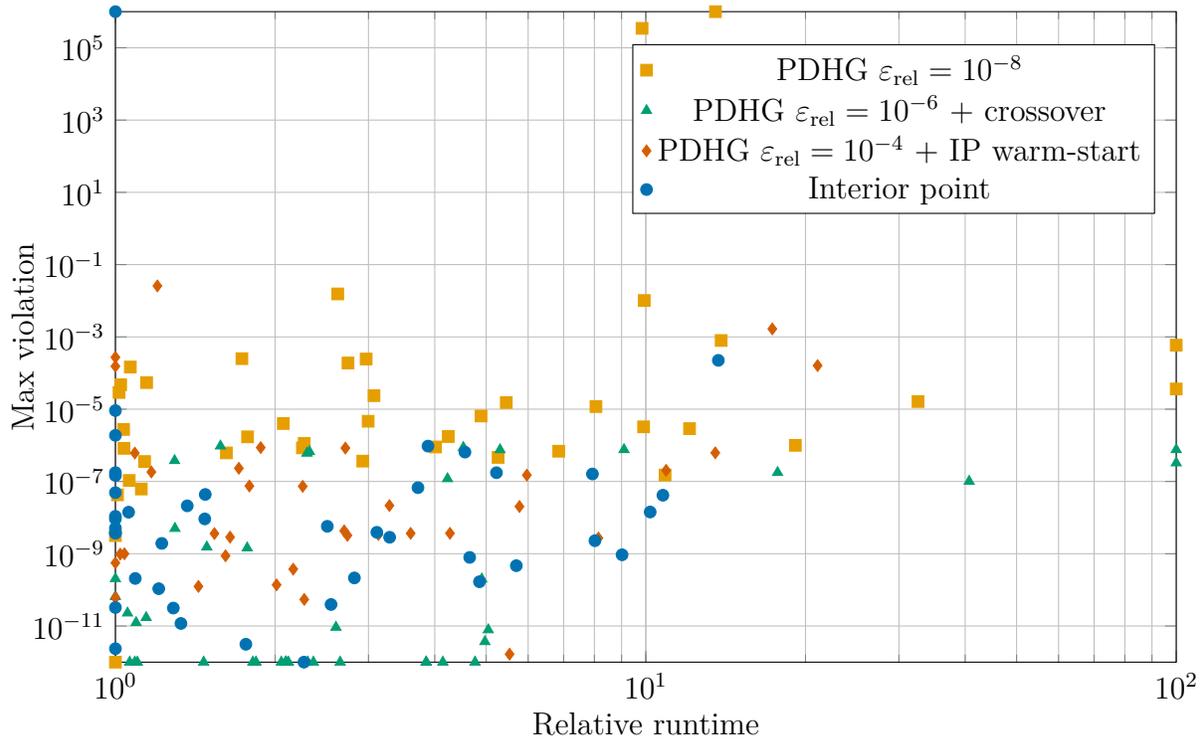
\begin{figure}[htbp]
\centering
\begin{tikzpicture}
\begin{loglogaxis}[
    width=0.95\linewidth,
    height=0.62\linewidth,
    grid=both,
    legend style={at={(0.98,0.95)},anchor=north east},
    mark size=2.2pt,
    xmin=1,
    xmax=1e2,
    ymin=1e-12,
    ymax=1e6,
    xlabel={Relative runtime},
    ylabel={Max violation},
    ]

\definecolor{OIblue}{RGB}{0,114,178}
\definecolor{OIorange}{RGB}{230,159,0}
\definecolor{OIteal}{RGB}{0,158,115}
\definecolor{OIvermillion}{RGB}{213,94,0}

\addplot[
    only marks,
    mark=square*,
    color=OIorange,
] table {
x y
1.016e+00 2.895e-05
1.038e+00 8.264e-07
1.062e+00 1.064e-07
4.890e+00 6.464e-06
2.927e+00 3.657e-07
2.073e+00 4.014e-06
1.775e+00 1.715e-06
1.387e+01 7.959e-04
1.135e+00 3.557e-07
1.023e+00 4.751e-05
1.732e+00 2.503e-04
1.000e+00 1.000e-12
1.353e+01 1.000e+06
1.009e+00 4.253e-08
2.252e+00 8.479e-07
1.067e+00 1.466e-04
2.745e+00 1.897e-04
1.145e+00 5.447e-05
1.000e+02 5.908e-04
1.000e+00 3.176e-09
2.625e+00 1.546e-02
5.457e+00 1.526e-05
4.242e+00 1.749e-06
5.267e+00 4.633e-07
1.914e+01 9.998e-07
2.970e+00 2.439e-04
3.260e+01 1.625e-05
8.053e+00 1.183e-05
9.932e+00 1.015e-02
1.119e+00 6.184e-08
4.013e+00 9.006e-07
2.996e+00 4.639e-06
6.852e+00 6.866e-07
1.000e+02 3.614e-05
3.073e+00 2.373e-05
1.087e+01 1.506e-07
2.270e+00 1.128e-06
1.209e+01 2.891e-06
1.036e+00 2.750e-06
1.619e+00 6.169e-07
9.904e+00 3.245e-06
9.838e+00 3.451e+05
};
\addlegendentry{PDHG $\epsrel = 10^{-8}$}

\addplot[
    only marks,
    mark=triangle*,
    color=OIteal,
] table {
x y
2.094e+00 1.000e-12
4.070e+01 1.010e-07
1.094e+00 1.243e-11
4.529e+00 8.706e-07
2.653e+00 1.000e-12
2.297e+00 6.046e-07
1.000e+00 2.044e-10
5.311e+00 7.550e-07
1.054e+00 2.285e-11
1.772e+00 1.445e-09
1.293e+00 3.781e-07
1.088e+00 1.000e-12
2.322e+00 6.688e-07
1.100e+00 1.000e-12
1.841e+00 1.000e-12
1.578e+00 9.553e-07
2.362e+00 1.000e-12
1.487e+00 1.534e-09
1.000e+02 7.585e-07
1.000e+02 3.256e-07
1.143e+00 1.715e-11
1.295e+00 5.006e-09
4.906e+00 2.001e-10
1.000e+00 6.544e-11
4.143e+00 1.000e-12
3.856e+00 1.000e-12
2.099e+00 1.000e-12
1.772e+01 1.770e-07
3.852e+00 1.000e-12
9.097e+00 7.588e-07
2.305e+00 1.000e-12
1.816e+00 1.000e-12
1.467e+00 1.000e-12
4.973e+00 3.688e-12
2.119e+00 1.000e-12
4.764e+00 1.000e-12
2.095e+00 1.000e-12
5.046e+00 7.850e-12
1.064e+00 1.000e-12
2.055e+00 1.000e-12
2.605e+00 9.218e-12
4.228e+00 1.192e-07
};
\addlegendentry{PDHG $\epsrel = 10^{-6}$ + crossover}

\addplot[
    only marks,
    mark=diamond*,
    color=OIvermillion,
] table {
x y
1.169e+00 1.830e-07
3.603e+00 3.658e-09
1.021e+00 9.774e-10
1.000e+00 5.583e-10
2.714e+00 8.398e-07
2.700e+00 4.291e-09
1.040e+00 9.968e-10
1.537e+00 3.617e-09
1.088e+00 6.114e-07
1.200e+00 2.575e-02
1.789e+00 7.491e-08
2.738e+00 3.220e-09
5.778e+00 2.028e-08
2.255e+00 7.326e-08
1.434e+00 1.246e-10
1.732e+01 1.664e-03
3.286e+00 2.151e-08
1.709e+00 2.316e-07
1.879e+00 8.658e-07
1.352e+01 6.197e-07
1.645e+00 2.882e-09
5.535e+00 1.662e-12
1.000e+00 2.727e-04
1.000e+00 6.025e-11
1.000e+00 1.539e-04
5.966e+00 1.500e-07
3.132e+00 3.408e-09
1.613e+00 8.762e-10
2.108e+01 1.610e-04
1.092e+01 2.031e-07
4.273e+00 3.672e-09
2.270e+00 5.441e-11
2.164e+00 3.743e-10
2.013e+00 1.377e-10
8.140e+00 2.722e-09
};
\addlegendentry{PDHG $\epsrel = 10^{-4}$ + IP warm-start}

\addplot[
    only marks,
    mark=*,
    color=OIblue,
] table {
x y
1.207e+00 1.077e-10
1.223e+00 1.933e-09
1.000e+00 5.103e-09
1.000e+00 3.817e-09
1.000e+00 3.690e-09
1.000e+00 4.904e-08
2.550e+00 3.945e-11
3.885e+00 9.512e-07
1.019e+01 1.428e-08
1.477e+00 4.366e-08
1.366e+00 2.116e-08
1.059e+00 1.412e-08
1.000e+00 1.000e+06
1.762e+00 3.126e-12
2.822e+00 2.139e-10
2.267e+00 1.000e-12
1.000e+00 1.084e-08
1.329e+00 1.177e-11
4.560e+00 6.482e-07
1.000e+00 9.153e-06
4.857e+00 1.678e-10
8.018e+00 2.296e-09
1.000e+00 1.758e-07
1.000e+00 9.226e-09
5.229e+00 1.744e-07
1.474e+00 9.202e-09
9.020e+00 9.356e-10
1.286e+00 3.138e-11
1.090e+00 2.070e-10
1.000e+00 1.460e-07
7.932e+00 1.610e-07
3.289e+00 2.856e-09
1.000e+00 2.359e-12
1.000e+00 1.885e-06
1.370e+01 2.244e-04
5.700e+00 4.669e-10
4.655e+00 7.922e-10
3.111e+00 3.915e-09
3.719e+00 6.723e-08
2.509e+00 5.741e-09
1.000e+00 3.248e-11
1.077e+01 4.160e-08
};
\addlegendentry{Interior point}
\end{loglogaxis}
\end{tikzpicture}
\caption{Runtimes and maximum violations for the models in the Mittelmann test set, for four solution approaches.  Runtime ratios larger than 100 are reported as 100.  Similarly, maximum violations smaller than $10^{-12}$ or larger than $10^{6}$ are reported as those values.}
\label{fig:count_vs_violation_mod}
\end{figure}
Figure~\ref{fig:count_vs_violation_mod} shows the same data, but at
the level of individual models.  Again, each dot represents one method
run on one model.  The x-axis shows runtime relative to the best time
achieved by any of the options (the baseline PDHG run is included
in this comparison, but not in the figure).  The y-axis shows maximum violation
for that model.

\subsection{Impact on Iteration Counts and Robustness}

We now expand our test set to get a better understanding of how robust
our PDHG/interior-point hybrid approach is, and how much it can be
expected to reduce the cost of the interior-point solve when compared
with a cold-start in a broader context.


\begin{table}[t]
\centering
\begin{tabular}{c|c|c|c|c}
\textbf{} & Models tried & Models solved & Success rate & IP iteration count ratio \\ \hline
$\epsrel = 10^{-4}$ & 43 & 35 & 81\% & 0.51 \\
$\epsrel = 10^{-6}$ & 42 & 37 & 88\% & 0.30 \\
\end{tabular}
\caption{Robustness of IP warm-start from different PDHG starting points on the Mittelmann LP test set (run on the Grace Hopper system).}
\label{tab:robustness_mittelmann}
\end{table}

\begin{table}[t]
\centering
\begin{tabular}{c|c|c|c|c}
\textbf{} & Models tried & Models solved & Success rate & IP iteration count ratio \\ \hline
$\epsrel = 10^{-4}$ & 2066 & 1725 & 83\% & 0.55 \\
$\epsrel = 10^{-6}$ & 1984 & 1726 & 87\% & 0.35 \\
\end{tabular}
\caption{Robustness of IP warm-start from different PDHG starting points on the broader test set (run on the AMD CPU cluster).}
\label{tab:robustness_broad}
\end{table}
Tables~\ref{tab:robustness_mittelmann} and~\ref{tab:robustness_broad}
show aggregated results over our two test sets.  We use PDHG to
compute a starting point to the tolerance indicated in the first
column; that solution is used to warm-start the interior-point solver.
The second column shows the number of instances from the test set
where the interior-point solver was invoked (it is not invoked if
presolve solves the problem or if PDHG is unable to find a solution
that meets the specified tolerances within the time limit).  The next
column shows the number of those models where the warm-started
interior-point solver was able to find a solution that met its
convergence criteria.  The third column shows the ratio between the
first two columns.
Finally, the last column shows the geometric mean of the ratio of
interior-point iterations for the hybrid approach against iterations
for a cold-started interior-point solve, captured on models that both
were able to solve.

As can be seen in these tables, results on the broader set are actually
quite similar to those on the smaller Mittelmann set.  The
warm-started interior-point solver in the hybrid approach is able to
solve around 85\% of the problems presented to it.  Interior-point
iteration counts are reduced by nearly a factor of two
with $\epsrel = 10^{-4}$ and nearly a factor of
three with $\epsrel = 10^{-6}$.

Table~\ref{tab:robustness_10iter} adds some additional color to these
results.  The first two rows repeat data from the previous tables:
they show the number of models presented to the interior-point solver
and the number of those that it was able to solve.  The third row is
new: it shows the number of models where the interior-point solver was
able to converge to desired tolerances in 10 or fewer iterations.  We
view this threshold as a strong indication that the starting point is
helping convergence quite a bit on that model.
\begin{table}[t]
\centering
\begin{tabular}{c|c|c|c}
  \textbf{} &
   \shortstack{\textbf{Interior point}\\ \textbf{cold-start}} &
   \shortstack{\textbf{PDHG}\\ $\epsrel = 10^{-4}$\\ \textbf{+ IP warm-start}} &
   \shortstack{\textbf{PDHG}\\ $\epsrel = 10^{-6}$\\ \textbf{+ IP warm-start}}\\ \hline
Number of models tried & 2179 & 2066 & 1984 \\
Number of models solved & 2090 & 1725 & 1726 \\
Number solved in $\leq 10$ iters  & 318 & 793 & 1029 \\
\end{tabular}
\caption{Interior-point iteration counts for different starting points.}
\label{tab:robustness_10iter}
\end{table}
The warm-started interior-point solver converges in 10 or fewer
iterations for roughly half of the models it solves, versus less
than 17\% for the cold-started solver.

While the results presented here hopefully demonstrate the promise of
our proposed hybrid approach, they also bring up the broader question
of how to fairly evaluate the performance of PDHG.  Our approach finds
a place on the efficient frontier in the performance-accuracy tradeoff
among the methods considered in Table~\ref{tab:accuracy_options} for
improving the accuracy of a PDHG solution, but a cold-started
interior-point method produced better results than any of these
approaches.  Our experience is that, despite significant anecdotal
wins, PDHG does not show significant advantages versus existing
alternatives when considered over a broad test set, even with very
loose convergence tolerances.

If PDHG wins are currently anecdotal, then perhaps we should look for
anecdotal wins for our proposed approach as well.  When comparing the
cold-started interior-point solver against our PDHG hybrid on the
Mittelmann test set, despite the significant loss over the entire set,
we found that the latter approach was more than 10\% faster in 6 cases
(and the former was more than 10\% faster in 16 cases).

\section{Discussion}

The hybrid approach proposed here is not the only option for
hybridizing PDHG and interior-point methods.  We also considered using
an abbreviated interior-point solve to compute a start point for PDHG.
Specifically, we used an iterative linear solver to compute direction
vectors in the early iterations of an interior-point solver (which is
both quick and GPU-friendly), terminating once a fixed total number of
iterative solver iterations were performed (either 1,000 or 10,000 in
our tests).  This approach provided PDHG a starting point with smaller
residuals and better centering than the trivial start that is
typically used.  We should note that while centering is crucial for
interior-point methods, it plays no role that we know of in PDHG
convergence arguments.  Hoping for a practical benefit was perhaps
wishful thinking, and indeed our experiments did not find a consistent
advantage.

We only used the GPU for PDHG in our tests.  Recent
work~\cite{cuopt25} has looked at using GPUs for interior-point
iterations as well.  While the advantages over a CPU appear to be
modest at this time, GPU systems are a fast-moving target, so results
could look quite different in the future.  For now, we simply note
that our hybrid approach allows you to obtain highly accurate
solutions using methods that can run (almost) entirely on a GPU.

\section{Conclusions}

This paper proposed a hybrid approach that combines a PDHG step to
compute a low-accuracy solution and a warm-started interior-point step
to significantly improve the accuracy of that solution.  While
warm-starting interior-point methods is a notoriously difficult
problem, we found that in this context a relatively simple approach
led to fast and reliable convergence.  The interior-point
solver required fewer than 10 iterations to converge for more than
half of the models in a broad test set.

One question that has proven a bit tough to answer is how to evaluate
the breadth of applicability of our approach.  For a specific model to
benefit, it would need the following to be true: (i) PDHG provides a
substantial performance advantage over the alternatives, (ii) the
application that consumes the solution requires more accuracy than
PDHG can provide (in the available time), and (iii) performing a few
interior-point iterations is not that expensive.  We have provided
data suggesting that the proposed method would provide significant
advantages in such situations.  The question of how often these
circumstances arise remains.

\section*{Acknowledgements}

We would like to thank Haihao Lu and David Torres Sanchez for their
helpful comments.

\bibliographystyle{plainnat}
\bibliography{references}

\end{document}